\documentclass[a4paper,11pt]{amsart}

\usepackage{amssymb,amsmath,epsfig,graphics,psfrag}
\usepackage[english,francais]{babel}

\usepackage[T1]{fontenc}

\newtheorem{theoreme}{Théorème}[section]
\newtheorem{proposition}{Proposition}[section]
\newtheorem{corollaire}{Corollaire}[section]
\newtheorem{lemme}[theoreme]{Lemme}

\theoremstyle{definition}
\newtheorem{definition}[theoreme]{Définition}

\theoremstyle{definition}
\newtheorem{remarque}[theoreme]{Remarque}

\numberwithin{equation}{section}

\newtheorem{theorem}{Theorem}
\newtheorem{lemma}{Lemma}
\newtheorem{e.proposition}{Proposition}
\newtheorem{e-definition}[theoreme]{Definition}

  \def\RR{{\mathbf{R}}}
 \def\TT{{\mathbf{T}}} 
  
 \def\ZZ{{\mathbf{Z}}}

\def\cA{{\mathcal{A}}}  
 \def\cE{{\mathcal{E}}}

 \def\cT{{\mathcal{T}}}

\def\fv{{\mathfrak{v}}}\def\fw{{\mathfrak{w}}}

\begin{document}
\sloppy

\title[Champs de vecteurs Morse-Smale sans singularité]{Classes d'homotopie de champs de vecteurs Morse-Smale sans singularité sur les fibrés de Seifert}

\author{Emmanuel Dufraine}
\address{Mathematics Institute, University of Warwick, CV4 7AL Coventry, U.K.}
%\curraddr{Department of Mathematics and Statistics, Case Western Reserve University, Cleveland, Ohio 43403}
\email{dufraine@maths.warwick.ac.uk}
%\thanks{\today}
\thanks{This research was supported by a Marie Curie Fellowship. (\today)}

%\author{Author Two}
%\address{Mathematical Research Section, School of Mathematical Sciences,Australian National University, Canberra ACT 2601, Australia}

%    General info

\selectlanguage{francais}
\begin{abstract} 
Nous considérons les applications d'une variété de dimension trois, compacte, orientable et sans bord dans la sphère $S^2$.
Nous donnons un critère permettant de décider si deux applications données sont homotopes, en fonction de l'ensemble des points où les applications sont égales et celui où elles sont opposées. Nous étendons ces résultats aux champs de vecteurs non-singuliers et aux champs de plans co-orientés sur les variétés de dimension trois.
Finalement, nous appliquons ce critère à l'étude des champs Morse-Smale non-singuliers sur les variétés de Seifert.
\end{abstract}

\maketitle

\vspace{-.5cm}
\tableofcontents

\selectlanguage{english}
\section*{English summary}
The first result of this article is a criterion for the existence of a homotopy between two smooth maps from a closed orientable 3-manifold to the 2-sphere. Let $f\colon M\to S^2$ be such a map, we put $c=[f^{-1}(y)]_{H_1(M,\ZZ)}$ (for $y$ a regular value) and we recall that the {\em maximal (free) divisor} of $c$ is 0 if $c$ is a torsion element, otherwise it is the largest non zero integer $p$ such that $c=pg$ with $g\in H_1(M,\ZZ)$.

If $g$ is an other map from $M$ to $S^2$, we pose $$C_+(f,g)=\{x\in M,\, f(x)=g(x)\}\mbox{ and}$$ $$C_-(f,g)=\{x\in M,\, f(x)=-g(x)\}.$$

Up to a small perturbation of $f$ and $g$, $C_+$ and $C_-$ are oriented embedded links in $M$.

\begin{lemma}\label{e.l.obst1geom}
For $f$ and $g$ two maps from $M$ to $S^2$, we have $$[C_+(f,g)]_{H_1 (M)}=c_f-c_g=[C_-(f,g)]_{H_1 (M)}\mbox{ (up to sign)}.$$ 
\end{lemma}

So, if $[C_-(f,g)]_{H_1 (M)}=0$, we also have $[C_+(f,g)]_{H_1 (M)}=0$ and the linking number between $C_+$ and $C_-$ is well-defined. Our criterion is given by the next proposition.

\begin{e.proposition}\label{e.p.ns}
Two maps $f$ and $g$ from $M$ to $S^2$ are homotopic if and only if
\begin{itemize} 
\item $[C_-(f,g)]_{H_1 (M)}=0$, which implies $c_f=c_g$\,; and, with $p$ the maximal divisor of $c_f$ and $c_g$,
\item $Link (C_+(f,g),C_-(f,g))=0$ modulo $2p$.
\end{itemize}
\end{e.proposition}

One of the main interests of this result is that it extends directly to non-singular vector fields on $M$. If $X$ and $Y$ are two non-singular vector fields on $M$, we take $$C_+(X,Y)=\{x\in M,\, X(x)=\lambda Y(x),\,\lambda> 0\}\mbox{ and}$$ $$C_-(X,Y)=\{x\in M,\, X(x)=\lambda Y(x),\,\lambda< 0\}.$$

As before, up to a small perturbation, $C_+$ and $C_-$ are oriented embedded links in $M$.
We denote by $\cE(X)\in H_1(M,\ZZ)$ the Poincaré dual of the Euler class of $X^\bot$. Proposition~\ref{e.p.ns} gives~:

\begin{e.proposition}\label{e.p.champ.ns}
Two non-singular vector fields $X$ and $Y$ on $M$ are homotopic if and only if
\begin{itemize} 
\item $[C_-(X,Y)]_{H_1 (M)}=0$, which implies $\cE(X)=\cE(Y)$\,; and, with $p$ the maximal divisor of $\cE(X)$ and $\cE(Y)$,
\item $Link(C_+(X,Y),C_-(X,Y))=0$ modulo $p$.
\end{itemize}
\end{e.proposition}

We apply this result to the study of non-singular Morse-Smale vector fields, giving a new proof of the following theorem which is a combination of results of Yano and Wilson~\cite{ya0,ya,wi}.

\begin{theorem}\label{e.t.nms}
For every Seifert manifold $M^3$, there exists an integer $n(M)$ such that every vector field on $M$ is homotopic to a non-singular Morse-Smale vector field with less than $n(M)$ periodic orbits.
\end{theorem}

\selectlanguage{francais}
\section{Introduction}
De nombreuses recherches actuelles portent sur l'étude de champs de 2-plans ou de champs de vecteurs partout non nuls, tangents à une variété de dimension trois (structures de contact, feuilletacts, champs de Morse-Smale non-singuliers~\dots~voir en particulier~\cite{cogiho,elth,hon,ku,mk}). La classification de certaines structures particulières (structures de contact tendues par exemple) dans les différentes classes d'homotopie est un sujet très actif actuellement. Le fibré tangent d'une variété de dimension trois, compacte, orientable est trivialisable\,; si la trivialisation est fixée, un champ de 2-plans co-orienté ou un champ de vecteurs non-singulier est alors uniquement associé à une application de $M$ dans la sphère $S^2$. 

Nous donnons dans cette note un critère géométrique pour décider si deux applications de $M$ dans $S^2$ sont homotopes ou non (Proposition~\ref{p.ns}). L'avantage de ce critère est qu'il se généralise aux champs de vecteurs non-singuliers (ou champs de plans co-orientés) et qu'il ne dépend pas du choix d'une trivialisation du fibré tangent de la variété.

Nous utilisons ce critère pour redémontrer de manière élémentaire (sans utiliser les décompositions en anses rondes d'Asimov et Morgan) le résultat suivant, conséquence des travaux de Yano~\cite{ya0,ya} et de Wilson~\cite{wi}~:

\begin{theoreme}\label{t.nms}
Pour chaque variété de Seifert $M^3$, il existe un nombre $n(M)$ tel que tout champ de vecteurs non-singulier sur $M$ est homotope à un champ de Morse-Smale non-singulier ayant au plus $n(M)$ orbites périodiques.
\end{theoreme}

Les champs de vecteurs Morse-Smale non-singuliers sur les variétés de dimension trois admettent beaucoup de propriétés intéressantes, malgré leur apparente simplicité. On pourra en particulier se reporter aux articles de Franks et Wada~\cite{fr,wa} concernant la topologie des orbites périodiques des champs de la sphère $S^3$. Les bifurcations des entrelacs d'orbites périodiques ont été étudiées dans~\cite{camavi}. Enfin des liens étonnants ont été mis en évidence avec les hamiltoniens intégrables dans~\cite{camanu}. 

Le Théorème~\ref{t.nms} montre en particulier que l'on peut toujours <<simplifier>> par une homotopie n'importe quelle dynamique sur une variété de Seifert. On voit aussi que le nombre d'orbites périodiques d'un champ de Morse-Smale non-singulier, que l'on peut interpréter comme une mesure de complexité pour ces champs, n'est pas relié à sa classe d'homotopie.

Dans~\cite[\S 6]{mk}, MacKay propose d'étudier l'influence de la géométrie de la variété (au sens de Thurston) sur l'existence de {\em dynamique compliquée}, à homotopie près (il emploie le terme {\em isotopy} pour l'homotopie de champs non-singuliers). En particulier, il rappelle que l'on peut rendre périodique, par une homotopie, le flot géodésique sur une surface. Le Théorème~\ref{t.nms} entra\^\i ne l'existence d'une homotopie du flot géodésique vers un Morse-Smale non-singulier.

Le critère d'homotopie et le processus de construction de champs de Morse-Smale présentés ici sont des généralisation de~\cite{du_qtds} où nous menions cette étude sur la sphère $S^3$. Nous espérons que les techniques employées ici pourront \^etre utilisées dans d'autres situations.

\subsection{Critère d'homotopie pour les applications}

On considère $M$ une variété de dimension trois, compacte, orientable, sans bord et nous nous intéressons aux applications lisses de $M$ dans la sphère $S^2$, à homotopie près. 

Si la variété de départ est la sphère $S^3$ (ou une sphère d'homologie), Hopf associe à une application $f$ un nombre, $H(f)$, en calculant l'enlacement entre deux images réciproques de valeurs régulières de $f$. Ce nombre ne dépend pas du choix des valeurs régulières et est invariant à homotopie de $f$ près. De plus, il classifie ces applications à homotopie près (c.f~\cite{mi} par exemple). 

Plus généralement, si $f$ est une application de $M$ dans $S^2$, l'image réciproque d'une valeur régulière $y$ est une sous-variété orientée de codimension 2 de $M$. Le choix d'une base du plan tangent à $S^2$ au point $y$ permet de trivialiser le fibré normal de $f^{-1}(y)$, on dit que cette sous-variété associée à $f$ est <<{\em framée}>>. Dans~\cite{po,po2} (voir aussi~\cite[\S 7]{mi}), Pontryagin montre que la variété framée associée à $f$ ne dépend pas, à {\em cobordisme framé} près, du choix de la valeur régulière ou du choix de la base du plan tangent à $S^2$. Sa classe de cobordisme framé est indépendante de $f$ dans sa classe d'homotopie et les classes de cobordisme framé des entrelacs framés d'une variété de dimension trois sont en bijection, par cette association, avec les classes d'homotopie des applications dans $S^2$. 

L'existence d'une homologie entre deux entrelacs d'une variété de dimension trois étant équivalente à l'existence d'un cobordisme entre ces entrelacs, on associe donc à une application $f$ une {\em classe caractéristique}, $c_f$, dans $H_1(M,\ZZ)$ qui est la classe d'homologie de $f^{-1}(y)$ pour $y$ une valeur régulière. Cette classe caractéristique ne dépend pas du choix de $y$ et est invariante si on change $f$ par une homotopie. D'après le résultat de Pontryagin, il reste donc à comprendre la partie <<framing>> pour caractériser la classe d'homotopie de $f$. Pour cela, nous rappelons qu'un élément $\tau$ de $H_1(M,\ZZ)$ est {\em de torsion} s'il existe un entier $k$ non nul tel que $k\tau=0$.

\begin{definition}
Le {\em diviseur maximal (libre)} d'une classe d'homologie $c$ de $H_1(M,\ZZ)$ est 
\begin{itemize}
\item le plus grand entier non nul $p$ vérifiant $c=p g$ pour $g$ dans $H_1(M,\ZZ)$ si $c$ n'est pas de torsion,
\item nul si $c$ est de torsion.
\end{itemize}

Un élément dont le diviseur maximal est égal à 1 est appelé {\em générateur} ou {\em élément primitif}.
\end{definition}

Pontryagin montre dans~\cite[\S4]{po} la proposition suivante (voir aussi~\cite[Theorem 6.2.7]{bepe}, \cite[Proposition 4.1]{go} et \cite[Proposition 2.1]{ku1} pour des preuves plus modernes). 

\begin{proposition}
La classe d'homotopie d'une application $f$ de $M$ dans $S^2$ est décrite par sa classe caractéristique $c_f$ dans $H_1(M,\ZZ)$ et un <<{\em degré de Hopf}>> $d$ dans un espace affine $\ZZ_{2p_f}$ où $p_f$ est le diviseur maximal de $c_f$.
\end{proposition}

Nous expliquons plus précisément cette description à la section \ref{s.hopf}.

Le problème de cette description est qu'il n'y a pas de manière canonique d'identifier l'espace affine $\ZZ_{2p_f}$ avec $\ZZ_{2p_f}$ (voir l'exemple au \S\ref{s.ex}). On peut montrer qu'une telle identification existe si $c_f$ est nulle ou de torsion\,; en particulier, on retrouve l'invariant de Hopf dans $\ZZ$ si $M=S^3$ (c.f.~\cite{po}).

Dans~\cite{ku1}, Kuperberg remarque après la preuve de la Proposition 2.2, que si l'on compare deux applications $f$ et $g$ ayant m\^eme classe caractéristique $c$, la différence des degrés de Hopf de $f$ et $g$ est un élément de $\ZZ_{2p_f}$ de manière canonique.

Le premier objectif de cette note est de donner une preuve géométrique de ce fait. Pour cela, nous définissons les ensembles $$C_+(f,g)=\{x\in M,\, f(x)=g(x)\}\mbox{ et}$$ $$C_-(f,g)=\{x\in M,\, f(x)=-g(x)\}.$$

On montre que $C_+$ et $C_-$ sont, quitte à perturber $f$ et $g$, des entrelacs orientés de $M$. Nous montrons alors~:

\begin{lemme}\label{l.obst1geom}
Pour $f$ et $g$ de $M$ dans $S^2$, on a $$[C_+(f,g)]_{H_1 (M)}=c_f-c_g=[C_-(f,g)]_{H_1 (M)}\mbox{ (au signe près)}.$$ 
\end{lemme}

L'ambiguité du signe est levée par le choix d'une orientation de $S^2$ et de $M$.

D'après le lemme précédent, si la classe d'homologie de $C_-(f,g)$ est nulle, celle de $C_+(f,g)$ l'est aussi; l'enlacement entre $C_+$ et $C_-$ est donc bien défini dans ce cas. On obtient alors le premier résultat de cette note~:
\begin{proposition}\label{p.ns}
Deux applications $f$ et $g$ de $M$ dans $S^2$ sont homotopes si et seulement si
\begin{itemize} 
\item $[C_-(f,g)]_{H_1 (M)}=0$, ce qui entra\^\i ne $c_f=c_g$\,; et, en notant $p$ le diviseur maximal de $c_f$ et $c_g$,
\item $Enl(C_+(f,g),C_-(f,g))=0$ modulo $2p$.
\end{itemize}
\end{proposition}

Si le diviseur $p$ est inconnu, on a la condition suffisante~:

\begin{corollaire}\label{c.s}
\'Etant données deux applications $f$ et $g$ de $M$ dans $S^2$, 
si $[C_-(f,g)]_{H_1 (M)}=0$ et $Enl(C_+(f,g),C_-(f,g))=0$, $f$ et $g$ sont homotopes.
\end{corollaire}

\begin{remarque}
On peut définir une {\em distance en classes d'homotopies} entre deux applications~:
\begin{itemize}
\item la classe d'homologie de $C_-$ donne la différence en classes d'homologie entre les classes caractéristiques\,;
\item Si la distance entre les classes caractéristiques est nulle, l'enlacement $Enl(C_+,C_-)$ donne la distance, dans $\ZZ_{2p}$, en classes d'homotopies ayant la m\^eme classe caractéristique entre les deux applications.  
\end{itemize}
\end{remarque}

\subsection{Critère d'homotopie pour les champs de vecteurs}

Si $X$ et $Y$ sont deux champs de vecteurs non-singuliers sur $M$, on pose $$C_+(X,Y)=\{x\in M,\, X(x)=\lambda Y(x),\,\lambda> 0\}\mbox{ et}$$ $$C_-(X,Y)=\{x\in M,\, X(x)=\lambda Y(x),\,\lambda< 0\}.$$

Comme pour les applications, le théorème de transversalité permet de supposer (après une éventuelle petite perturbation de $X$ et $Y$) que $C_+$ et $C_-$ sont des entrelacs orientés plongés dans $M$.

On note $\cE(X)\in H_1(M,\ZZ)$ le dual de Poincaré de la classe d'Euler de $X^\bot$. Pour les champs de vecteurs, la Proposition~\ref{p.ns} donne le résultat~:

\begin{proposition}\label{p.champ.ns}
Deux champs de vecteurs non-singuliers $X$ et $Y$ sur $M$ sont homotopes si et seulement si
\begin{itemize} 
\item $[C_-(X,Y)]_{H_1 (M)}=0$, ce qui entra\^\i ne $\cE(X)=\cE(Y)$\,; et, en notant $p$ le diviseur maximal de $\cE(X)$ et $\cE(Y)$,
\item $Enl(C_+(X,Y),C_-(X,Y))=0$ modulo $p$.
\end{itemize}
\end{proposition}

La Proposition~\ref{p.champ.ns} est une généralisation de~\cite[Proposition 1.1]{neru} (voir aussi \cite[Lemma 23]{du_qtds}) où ce résultat est prouvé dans le cas de la sphère $S^3$.

\begin{remarque}
Pour les champs de plans co-orientés (feuilletages, structures de contact, feuilletacts, \dots), il suffit de remplacer $C_+$ (resp. $C_-$) par l'ensemble des points où les deux champs de plans co\"\i ncident avec la m\^eme orientation (resp. orientation opposée) pour que la proposition précédente soit valide.
\end{remarque}

\subsection{Champs de vecteurs Morse-Smale non-singuliers}

\begin{definition}
Un champ de vecteurs $X$ est de type {\em Morse-Smale non-singulier} sur $M$ si
\begin{itemize}
\item aucun point de $M$ n'est fixé par le flot $\phi$ de $X$, 
\item l'ensemble non-errant de $\phi$ est réduit à un nombre fini d'orbites périodiques hyperboliques~et
\item les variétés stables et instables des orbites périodiques s'intersectent transversalement.
\end{itemize}
\end{definition}

Asimov montre dans~\cite{as1} que sur toute variété de dimension supérieure ou égale à quatre et de caractéristique d'Euler nulle, tout champ de vecteurs non-singulier est homotope à un champ de Morse-Smale non-singulier. Morgan montre ensuite, \cite{mo}, que ce résultat ne peut \^etre vrai en dimension trois puisque beaucoup de variétés n'admettent aucun champ de Morse-Smale non-singulier. En généralisant la notion de décomposition en anses rondes d'Asimov, il montre en effet que les seules variétés de dimension trois orientables, à bord torique et premières pour la décomposition en somme connexe admettant des champs Morse-Smale non-singuliers sont les variétés graphées (recollement de variétés de Seifert le long de leurs bords). 

Sur les variétés de dimension trois admettant des champs Morse-Smale non-singuliers, Yano donne dans~\cite{ya} une caractérisations des classes d'homotopies admettant des champs de Morse-Smale non-singuliers. En particulier, il montre que sur les variétés de Seifert, toute classe d'homotopie de champs de vecteurs admet un champ Morse-Smale non-singulier. Sur la sphère $S^3$, il existe un champ de vecteurs Morse-Smale non-singulier avec au plus 6 orbites périodiques dans chaque classe d'homotopie (voir~\cite{wi} et \cite{du_qtds}). Yano en déduit l'existence d'un nombre $n(M)$, tel qu'il existe un champ de type Morse-Smale avec au plus $n(M)$ orbites périodiques dans chaque classe d'homotopie de champs de vecteurs possible sur $M$.

L'idée de Yano est tout d'abord de montrer que s'il existe un champ de Morse-Smale ayant une classe d'Euler donnée, alors toutes les classes d'homotopie ayant cette classe d'Euler  admettent un champ de Morse-Smale. Puis il montre, sur les variétés de Seifert, que l'on peut construire un champ de Morse-Smale ayant n'importe quelle classe d'Euler. 

Nous utilisons la m\^eme stratégie, mais quand Yano utilise la somme connexe avec $S^3$ pour réaliser la première étape, nous utilisons le critère donné à la Proposition~\ref{p.champ.ns} pour construire les champs de vecteurs (Proposition~\ref{p.enlmorse}). En outre, nous n'utilisons pas la décomposition en anses rondes de la variété, élément essentiel de la preuve de Yano. Enfin, nous précisons les résultats de~\cite{ya0}, afin de prouver l'existence du nombre $n(M)$. 

L'avantage de notre approche est d'une part que nous construisons explicitement des champs de vecteurs Morse-Smale dans chaque classe d'homotopie\,; d'autre part, les orbites périodiques des champs de vecteurs que nous construisons ne sont pas homologues à zéro. En particulier, ces champs sont tous transverses à un feuilletage gr\^ace à un résultat de Goodman (\cite{go1}, voir aussi~\cite{ya2}). Cela nous amène à discuter d'une classification des entrelacs essentiels (sans composante homologue à zéro) similaire à celle de Wada~\cite{wa} dans le cas de la sphère (section~\ref{s.discussion}).

\medskip

{\bf Remerciements~:} Je tiens à remercier Pierre Derbez pour son aide <<homologique>> et l'intér\^et qu'il a porté à ce travail. Je remercie aussi Max Forester, Robert MacKay et Colin Rourke pour de nombreuses discussions très motivantes ainsi que Daniel Lines pour son aide précieuse.

\section{Entrelacs framés et modèles de Pontryagin}\label{s.hopf}

Nous considérons $M$ une variété de dimension trois, compacte et orientée.

\begin{definition}
Un entrelac $L$ dans $M$ est {\em framé} s'il est orienté et s'il existe une trivialisation $\fv$ du fibré normal à $L$, compatible avec l'orientation de $M$.

Deux entrelacs framés $(L,\fv)$ et $(\tilde L,\tilde \fv)$ sont {\em frame-cobordants} s'il existe une surface $S$, plongée dans $M\times I$, telle que $S$ rencontre le bord de $M\times I$ transversalement, $\partial S=L\times\{0\}\cup -\tilde L\times\{1\}$ et s'il existe $\mathfrak V$ une trivialisation du fibré normal de $S$ co\"\i ncidant avec $\fv$ et $-\tilde \fv$ le long de $L\times\{0\}$ et de $-\tilde L\times\{1\}$.
\end{definition}

\begin{remarque}
Dans la définition précédente, la surface $S$ est nécessairement orientable.
\end{remarque}

La donnée d'une trivialisation du fibré normal à un entrelac $L$ dans $M$ (resp. à une surface (orientable) $S$ dans $M\times I$) est équivalente à la donnée d'un champ de vecteurs non-singulier (tangent à $M$ (resp. $M\times I$)) normal à $L$ (resp. $S$).

Réciproquement, étant donnée une trivialisation $\fv$ du fibré normal de $L$ dans $M$ (resp. d'une surface $S$ dans $M\times I$), nous appellerons {\em champ de vecteurs constant dans $\fv$} un champ de vecteurs au voisinage de $L$ (resp. $S$) dont la restriction à $L$ (resp. $S$) est envoyée sur une constante par $\fv$.

\subsection{Enlacement d'entrelacs homologues à zéro}\label{s.enlace}
On considère $K$ et $L$ deux entrelacs disjoints et orientés, homologues à zéro dans $M$. Soit $\sigma_K$ une 2-cha\^\i ne bordant $K$ ($\sigma_K$ existe car $K$ est homologue à zéro), $\sigma_K$ est orientée de sorte que l'orientation qu'elle induit sur son bord soit celle de $K$. On compte le nombre algébrique de points d'intersection entre $\sigma_K$ et $L$ (génériquement, ils s'intersectent transversalement) et on pose $Enl(K,L)=\sigma_K \cdot L$. Ce {\em nombre d'enlacement} ne dépend pas du choix de $\sigma_K$ car $L$ est aussi homologue à zéro\,; on montre aussi que $Enl(K,L)=Enl(L,K)$.

Soit $K$ un n\oe ud orienté de $M^3$ et $\cT$ une paramétrisation de son voisinage tubulaire (i.e. un difféomorphisme $\cT\colon D^2\times S^1\to M$ tel que  $\cT(D^2\times S^1)$ soit un voisinage de $K=\cT(\{0\}\times S^1)$). L'image par $\cT$ d'un cercle $\{x\}\times S^1$ sera appelée un {\em 1-c\^able} de $K$. Plus généralement~:

\begin{definition}
Un $n$-c\^able $L$ (pour $n$ un entier non nul) de $K$ est un n\oe ud dans un voisinage tubulaire de $K$, tel que le nombre algébrique d'intersections de $L$ avec un disque transverse à $K$ soit égal à $n$.
\end{definition}

Si $K$ est homologue à zéro, on peut définir plus précisément, pour $p$ et $q$ deux entiers premiers entre eux ($p$ non nul) :
\begin{definition}
Un $(p,q)$-c\^able $L$ de $K$ est un $p$-c\^able de $L$ tel que l'enlacement entre $L$ et $K$ soit égal à $q$.
\end{definition}

L'image d'un n\oe ud $(p',q')$ de $\partial D^2 \times S^1$ par une paramétrisation du voisinage tubulaire de $K$ donne un $(p',q'+k)$-c\^able de $K$. Comme $K$ est homologue à zéro, il existe une paramétrisation du voisinage tubulaire de sorte que $k$ soit nul.

Si $\gamma$ est un entrelac dans $M$ et que $\fv$ est une trivialisation de son fibré normal, on note $\fv\gamma$ l'entrelac obtenu en poussant $\gamma$ le long d'un champ de vecteurs constant dans $\fv$.
 
\begin{definition}
Pour $\gamma$ un n\oe ud homologue à zéro dans $M$, on note $\fw_n$ une trivialisation de son fibré normal telle que $\fw_n\gamma$ soit un $(1,n)$-c\^able de $\gamma$.
\end{definition}

\subsection{Classification des entrelacs framés}

Soient $(L,\fv)$ et $(\tilde L,\tilde \fv)$ deux entrelacs framés, tels que $L$ et $\tilde L$ sont  homologues\,; on note $c$ la classe d'homologie de $L$ et de $\tilde L$ et $p$ son diviseur maximal.

Soit $S\subset M\times I$ un cobordisme réalisant cette homologie, on suppose pour l'instant que $S$ est connexe. Soit $D$ un disque de $S$, il est facile de voir que l'on peut trivialiser le fibré normal de $S\setminus D$ en étendant $\fv$ et $-\tilde \fv$, on note $\mathfrak V_{S\setminus D}$ une telle trivialisation. On considère un petit disque dans le complément de $L\cup\tilde L$ dans $M\times \{1\}$, on note $\gamma$ son bord.

\begin{figure}[htb]
\psfrag{L0}{$L\subset M\times \{0\}$}\psfrag{S}{$S$}\psfrag{D}{$D$}\psfrag{g}{$\gamma$}\psfrag{L1}{$\tilde L \subset M\times \{1\}$}
\centerline{\includegraphics{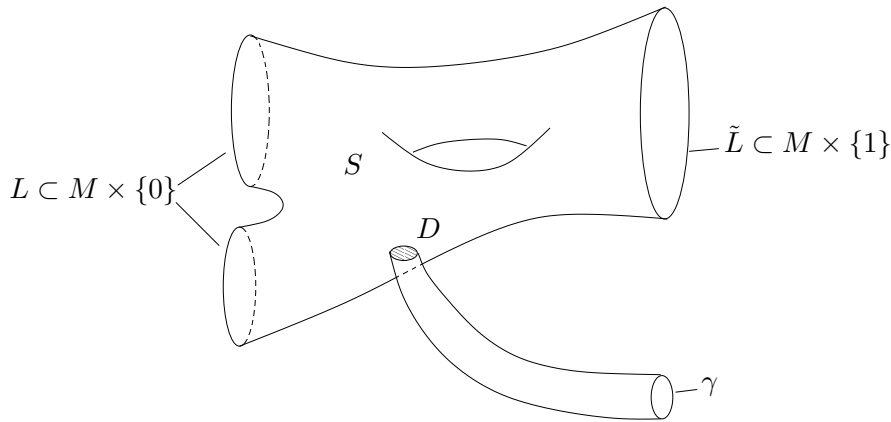}}
\caption{\label{f.cobord} La surface $S$ privée de $D$ et son cylindre}
\end{figure}

Il est facile de voir qu'il existe un entier $n$ tel que  $(L,\fv)$ et $(\tilde L\coprod \gamma,\tilde \fv\coprod \fw_{n})$ sont frame-cobordants, le cobordisme étant obtenu en collant un cylindre à $S$ le long de $\partial D$ et en prolongeant $\mathfrak V_{S\setminus D}$ à ce cylindre (cette construction peut \^etre trouvée aussi dans~\cite[\S 3.3.3]{ge}, voir la figure~\ref{f.cobord}).

Comme nous avons supposé $S$ connexe, ce nombre $n$ est indépendant du choix de $D$ dans $S$. De plus, on a le lemme suivant.

\begin{lemme}
Le nombre $n$ ne dépend pas du choix de $S$ ou de $\mathfrak V_{S\setminus D}$ modulo $2p$.
\end{lemme}

\begin{proof}
Notons $n(S,\mathfrak V)$ le nombre obtenu avec $S$ et $\mathfrak V_{S\setminus D}$ pour $D$ un disque de $S$. Considérons $S'$ un autre cobordisme entre $L$ et $\tilde L$, pour $D'$ un disque de $S'$ et $\mathfrak V'_{S'\setminus D'}$ une trivialisation du fibré normal de $S'\setminus D'$, nous notons $n(S',\mathfrak V')$ le nombre obtenu. Notre objectif est de prouver $n(S,\mathfrak V)\equiv n(S',\mathfrak V') \mod 2p$.

\begin{itemize}
\item Soit $\mathfrak V_0$ un champ de vecteurs, tangent à $M\times I$, définit sur $S\setminus D$, constant dans $\mathfrak V$. Le champ de vecteurs $\mathfrak V_0$ est partout normal à $S\setminus D$, on le prolonge en un champ $X_0$ au voisinage de $S\setminus D$ et on note $S_0$ la surface obtenue en poussant $S\setminus D$ le long de $X_0$\,; on note $\gamma_0$ la composante de bord de $S_0$ correspondant à $\partial D$. On remarque que les surfaces $S\setminus D$ et $S_0$ ne s'intersectent jamais.

\begin{lemme}  Il est possible de coller un disque à $S_0$, le long de $\gamma_0$, pour obtenir une surface $S_1$ ayant un nombre d'intersection avec $S$ égal à $n(S,\mathfrak V)$.
\end{lemme} 

\begin{proof}
Pour cela, il suffit de remarquer que, par définition de $n(S,\mathfrak V)$, on peut coller un cylindre $C$ à $S\setminus D$ le long de $\partial D$ de sorte que son autre bord $\gamma$ soit dans $M\times \{1\}$ le bord d'un disque $F$. Nous pouvons étendre $\mathfrak V_0$ le long de ce cylindre et ainsi, en poussant $\gamma$ le long de $\mathfrak V_0$, on obtient un n\oe ud s'enla\c cant $n(S,\mathfrak V)$ fois autour de $\gamma$. 

La surface $S$ est homotope à l'union $\overline S=(S\setminus D)\cup C \cup F$, de m\^eme, la surface $S_0$ est homotope à la surface $\overline {S_0}$ obtenue en poussant $(S\setminus D)\cup C$ le long de $\mathfrak V_0$. Les surfaces $\overline S$ et $\overline {S_0}$ se coupent exactement en $n(S,\mathfrak V)$ points, tous dans $F$. Il est maintenant facile de coller un disque à $\overline {S_0}$ ne coupant pas $\overline S$ ailleurs qu'en ces points. Nous obtenons ainsi la surface $S_1$, en faisant une homotopie dans l'autre sens pour que $S_1$ ne rencontre $M\times\{0,1\}$ qu'en son bord, le nombre d'intersection étant préservé par homotopie, on obtient le lemme. 
\end{proof}

\item Nous recollons deux copies de $M\times I$ le long de leurs bords (en renversant l'orientation d'une des copies) afin d'obtenir la variété de dimension 4, sans bord, $M\times S^1$. On recolle les surfaces $S$ et $-S'$ le long de $L\times\{0\}$ et $\tilde L\times\{1\}$ pour obtenir une surface fermée $\Sigma$ dans $M\times S^1$. De m\^eme, les surfaces $S_1$ et $-S_1'$ construitent comme ci-dessus se recollent le long des 1-c\^ables de $L\times\{0\}$ et $\tilde L\times\{1\}$ obtenus en poussant $L\times\{0\}$ et $\tilde L\times\{1\}$ par des champs de vecteurs constants dans les trivialisations $\fv$ et $\tilde \fv$ respectivement. On obtient ainsi une surface $\Sigma_1$. En prenant garde aux orientations, on voit que $\Sigma$ et $\Sigma_1$ s'intersectent en $n(S,\mathfrak V)-n(S',\mathfrak V')$ points (au signe près).

\item Nous terminons la preuve comme~\cite[Theorem 6.2.7]{bepe} ou \cite[Proposition 4.1]{go}. La classe d'homologie (dans $H_2(M\times S^1)$) de $\Sigma$ est la m\^eme que celle de $\Sigma_1$ et elle s'écrit, d'après la formule de K\"unneth~:
$$[\Sigma]=[a\times S^1]+[b\times\{w\}]$$ avec $a$ un 1-cycle de $M$ et $b$ un 2-cycle de $M$. On voit facilement que $[a]=c$ la classe d'homologie de $L$ et $\tilde L$.

Ainsi, $[\Sigma]\cdot [\Sigma_1]= [a\times S^1]\cdot[a\times S^1] + 2[a\times S^1]\cdot[b\times\{w\}] +[b\times\{w\}]\cdot[b\times\{w\}]$. D'où $$n(S,\mathfrak V)-n(S',\mathfrak V') = [\Sigma]\cdot [\Sigma_1] = 2c\cdot[b]$$ avec le dernier produit dans $H_\star(M)$. 

On en déduit que si $c$ est de torsion, $2c\cdot[b]$ est de torsion dans $H_3(M)=\ZZ$, ce qui entra\^\i ne $n(S,\mathfrak V)-n(S',\mathfrak V') = 0$.
Si $c$ n'est pas de torsion, $p$ divise $c$ et $n(S,\mathfrak V)-n(S',\mathfrak V')\equiv 0 \mod 2p$. \qedhere
\end{itemize}
\end{proof}

\begin{definition}\label{d.diff}
La {\em différence de degré de Hopf} entre $(L,\fv)$ et $(\tilde L,\tilde \fv)$ est l'élément de $\ZZ_{2p}$ défini par $n((L,\fv),(\tilde L,\tilde \fv))=n(S,\mathfrak V) \mod 2p$ pour $S$ et $\mathfrak V$ comme ci-dessus.
\end{definition}

Le théorème de Hopf sur les applications de $S^3$ dans $S^2$ peut se généraliser pour obtenir~:
\begin{lemme}
Deux n\oe uds $L$ et $\tilde L$ homologues à zéro dans $M$, dont le fibré normal est trivialisé par $\fw_n$ et $\fw_{\tilde n}$ respectivement, sont frame-cobordants si et seulement si $n=\tilde n$. De plus, on a la formule $n((L,\fw_n),(\tilde L,\fw_{\tilde n}))=n -\tilde n$, au signe près.
\end{lemme}

On en déduit les deux résultats suivants. Le premier nous permet de déduire que l'hypothèse <<$S$ connexe>>, dans la définition~\ref{d.diff} de la différence de degré de Hopf, ne restreint pas la portée de cette définition.
 
\begin{lemme}
Si deux entrelacs framés $(L,\fv)$ et $(\tilde L,\tilde \fv)$ sont frame-cobordants, on peut trivialiser le fibré normal de tout cobordisme $S$ orientable pour obtenir un cobordisme framé dès que $S$ est connexe.
\end{lemme}

\begin{lemme}\label{l.entfra}
Deux entrelacs framés $(L,\fv)$ et $(\tilde L,\tilde \fv)$ sont frame-cobordants si et seulement si
\begin{itemize}
\item $[L]=[\tilde L]$ et, en notant $p$ le diviseur maximal de $[L]=[\tilde L]$,
\item $n((L,\fv),(\tilde L,\tilde \fv))\equiv 0 \mod 2p$.
\end{itemize}
\end{lemme}

Pontryagin montre que les classes de cobordisme framé d'entrelacs framés de $M$ sont en bijection avec les classes d'homotopie des application de $M$ dans $S^2$. En conséquence du Lemme~\ref{l.entfra} nous obtenons (voir \S\ref{s.ex} pour un exemple)~:

\begin{theoreme}[\cite{po}, \S 4]
La classe d'homotopie de $f\colon M \to S^2$ est caractérisée par~:
\begin{itemize}
\item sa classe caractéristique $c_f=[f^{-1}(y)]_{H_1(M)}$ et
\item un <<degré de Hopf>> dans un $\ZZ_{2p_f}$ affine dans $\ZZ$.
\end{itemize}
\end{theoreme}

\subsection{Modèles de Pontryagin}
Soit $c$ un élément de $H_1 (M)$ et $\gamma$ un n\oe ud réalisant cette classe d'homologie. Nous considérons $T(\gamma)$ un voisinage tubulaire de $\gamma$, identifié avec $S^1\times D^2$ par un difféomorphisme $\cT$. 

Nous paramétrons le cercle $S^1$ par $\omega \in [0,2\pi]$, le disque $D^2$ avec les coordonnées polaires usuelles $(r,\theta)$ avec $r\in [0,1]$ et $\theta\in \RR/2\pi\ZZ$. On identifie la sphère $S^2$ privée du p\^ole Sud (noté $S$) avec le disque $D^2$, en envoyant le p\^ole Nord sur $0$ et les méridiens sur les rayons $\{(r,\theta_0),r \in [0,1[\}$. 

Pour $n\in \ZZ$, on définit l'application <<rajoutant $n$ twists à droite aux 1-c\^ables de $\gamma$ donnés par $\cT$>> par 
$$P_{c,n}\colon M \to S^2 \mbox{ par } \left\{\begin{array}{ll} 
P_{c,n} \equiv S&\mbox{sur } M\setminus \overline{T(\gamma)}\\
P_{c,n}(\cT^{-1}(w,r,\theta))=\cT(w,r,n\omega + \theta) &\mbox{sur } T(\gamma).
\end{array}
\right.$$

On vérifie facilement que $P_{c,n}$ est continue, les points de $S^2\setminus\{S\}$ sont des valeurs régulières et $P_{c,n}^{-1}(N)=\gamma$.

\begin{theoreme}[Construction de Pontryagin]\label{t.pont}
Toute application $f$ de $M$ dans $S^2$ telle que $c_f=c$ est homotope à un $P_{c,n}$.

De plus $P_{c,n}$ est homotope à $P_{c,n'}$ si et seulement si $n\equiv n' \mod 2p$, où $p$ est le diviseur maximal de $c$.
\end{theoreme}

\begin{remarque}
Dans la construction précédente, l'utilisation des p\^oles Nord et Sud n'est pas essentielle, on peut construire un modèle ayant les m\^emes propriétés avec n'importe quelle paire de points distincts de $S^2$. 
\end{remarque}

\section{Comparaison d'applications}
On définit $C_+(f,g)$ (resp. $C_-(f,g)$) comme l'ensemble des points de $M$ où $f=g$ (resp $f=-g$). Quitte à perturber légèrement $f$ et $g$ (ce qui ne change pas les applications à homotopie près), le théorème de transversalité nous permet de supposer que $C_+(f,g)$ et $C_-(f,g)$ sont des sous-variétés plongées dans $M$, de codimension 2. 

En choisissant une orientation sur $M$ et sur $S^2$, les entrelacs $C_+(f,g)$ et $C_-(f,g)$ héritent d'une orientation naturelle de la fa\c con suivante~: on associe à $f$ un plongement $F\colon M\hookrightarrow M\times S^2$ qui envoie le point $x$ sur le couple $(x,f(x))$. On note $\pi$ la projection de $M\times S^2$ sur le premier facteur $M$. Génériquement, les images $F(M)$ et $G(M)$ s'intersectent transversalement dans $M\times S^2$, le long d'une sous-variété de dimension 1, naturellement orientée, qui est en bijection via $\pi$ avec $C_+(f,g)$. On peut tenir le m\^eme raisonnement pour orienter $C_-(f,g)$ avec l'intersection des images $F(M)$ et $-G(M)=(M,-g(M))$.

La classe d'homologie (dans $H_1(M,\ZZ)$) de $C_+(f,g)$ et de $C_-(f,g)$ est invariante si on modifie $f$ ou $g$ par homotopie. D'autre part, $C_+(f,g)$ et $C_-(f,g)$ ne s'intersectent jamais, on en déduit que la classe d'homologie de $C_+(f,g)$ (resp.  de $C_-(f,g)$) dans $H_1(M\setminus C_-(f,g),\ZZ)$ (resp. $H_1(M\setminus C_+(f,g),\ZZ)$) est aussi invariante par homotopie des applications $f$ et $g$.

\subsection{Première obstruction}
Si $f$ et $g$ sont homotopes, il existe une application $\tilde f$, homotope à $f$, telle que $C_-(\tilde f,g)=\emptyset$\,; ainsi $[C_-(f,g)]_{H_1 (M)}=0$ est une condition nécessaire pour que $f$ et $g$ soient homotopes.

\subsubsection{Interprétation géométrique}
L'interprétation de cette première obstruction est donnée par le Lemme~\ref{l.obst1geom} que nous rappelons ici.

\begin{lemme}\label{l.obst1geom2}
Pour $f$ et $g$ de $M$ dans $S^2$, on a $$[C_-(f,g)]_{H_1 (M)}=c_f - c_g =[C_+(f,g)]_{H_1 (M)}\mbox{ (au signe près)}.$$ 
\end{lemme}

\begin{proof}
On considère $y \in S^2$ une valeur régulière commune à $f$, $g$ et $-g$. Comme $y$ est valeur régulière de $g$, $-y$ est valeur régulière de $-g$\,; on note $\gamma_{f}=f^{-1}(y)$, $\gamma_{-g}=(-g)^{-1}(y)$ et $\beta_{-g}=(-g)^{-1}(-y)$, ils sont orientés comme images réciproques de $y$ et $-y$ par $f$ et $-g$ respectivement. Remarquons que la classe d'homologie de $\gamma_{-g}$ et de $\beta_{-g}$ est l'opposée de $c_g$~: $[\gamma_{-g}]=[\beta_{-g}]=c_{-g}=-c_g$. 

On peut supposer, en faisant une homotopie sur $f$ et $g$, que $\gamma_f$ et $\beta_{-g}$ sont disjoints et que $f\equiv -y$ sur le complément d'un voisinage tubulaire de $\gamma_f$ et $-g\equiv y$ sur le complément d'un voisinage tubulaire de $\beta_{-g}$. Ceci entra\^\i ne que $C_-(f,g)=\gamma_f\cup\beta_{-g}$ et que $[C_-(f,g)]_{H_1 (M)}=c_f + c_{-g}=c_f - c_g$.

De m\^eme, on peut modifier $f$ et $g$ de sorte que cette fois $f\equiv -y$ sur le complément d'un voisinage tubulaire de $\gamma_f$ et $-g\equiv -y$ sur le complément d'un voisinage tubulaire de $\gamma_{-g}$. On obtient $C_+(f,g)=\gamma_f\cup\gamma_{-g}$, ainsi $[C_+(f,g)]_{H_1 (M)}=c_f + c_{-g}=c_f - c_{g}$.

L'ambiguité sur le signe provient du choix de l'orientation de $M$ et de $S^2$, un choix différent change l'orientation globale de $C_+$ et $C_-$\,; par suite, il change le signe de leur classe d'homologie. 
\end{proof}

\subsection{Seconde obstruction}

Si $M$ est la sphère $S^3$, la première obstruction est toujours nulle, on montre facilement (en adaptant les résultats de~\cite[Proposition 1.1]{neru} et \cite[Lemma 23]{du_qtds}) que deux applications $f$ et $g$ sont homotopes si et seulement si l'enlacement entre $C_+(f,g)$ et $C_-(f,g)$ est nul. 

Supposons maintenant que la première obstruction à l'existence d'une homotopie entre $f$ et $g$ est nulle ($[C_-(f,g)]_{H_1 (M)}=0$) et de plus, supposons que l'on connait $p$, le diviseur maximal de $c_f = c_{g}$. En rappelant que $[C_+(f,g)]_{H_1 (M)}=[C_-(f,g)]_{H_1 (M)}$ (Lemme~\ref{l.obst1geom2}), nous montrons que dans le cas général~:
\begin{lemme}\label{l.obst2}
Les applications $f$ et $g$ comme ci-dessus sont homotopes si et seulement si l'enlacement entre $C_+(f,g)$ et $C_-(f,g)$ est nul modulo $2p$.
\end{lemme}
\begin{proof}
Considérons $f$ et $g$ telles que $c_f = c_g$\,; soient $\gamma_1$ et $\gamma_2$ deux n\oe uds disjoints de $M$ dont l'homologie est $c_f = c_g$ et tels qu'il existe un anneau  $\cA$ plongé dans $M$, tel que $\partial A=\gamma_1\cup \gamma_2$ (sans tenir compte de l'orientation). On considère $U_1$, $U_2$ des voisinages tubulaires de $\gamma_1$ et $\gamma_2$ respectivement ($U_1$ et $U_2$ sont disjoints), paramétrés de sorte que la trace de $\cA$ dans $U_1$ (resp. $U_2$) soit l'image de $J\times S^1$ avec $J$ un segment de $D^2$. En notant $\star$ un point de $S^2$, distinct des deux p\^oles, la section précédente nous permet de construire deux modèles de Pontryagin $P_{k_1}$ et $P_{k_2}$ tels que $P_{k_1}^{-1}(N)=\gamma_1$, $P_{k_1}|_{M\setminus\overline{U_1}}\equiv S$, $P_{k_2}^{-1}(\star)=\gamma_2$, $P_{k_2}|_{M\setminus\overline{U_2}}\equiv -\star$ et $P_{k_1}$ est homotope à $f$, $P_{k_2}$ est homotope à $g$. Nous remarquons que pour tout point $y$ de $S^2$, différent de $N$ et $S$, $P^{-1}_{k_1}(y)$ rencontre $\cA$ en $k_1$ points. Et pour tout point $y$ de $S^2$, différent de $\star$ et $-\star$, $P^{-1}_{k_2}(y)$ rencontre $\cA$ en $k_2$ points.

Comme $U_1$ et $U_2$ sont disjoints, on montre facilement que $$C_+(P_{k_1},P_{k_2})=P_{k_1}^{-1}(\star)\cup P_{k_2}^{-1}(N)\mbox{ et}$$  $$C_-(P_{k_1},P_{k_2})=P_{k_1}^{-1}(-\star)\cup P_{k_2}^{-1}(S).$$ 

De plus, comme $c_f = c_g$, le Lemme~\ref{l.obst1geom} implique que $[C_+]=[C_-]=0$. D'autre part, $P_{k_1}^{-1}(\star)$, $P_{k_2}^{-1}(N)$, $P_{k_1}^{-1}(-\star)$ et $P_{k_2}^{-1}(S)$, orientés comme images réciproques de points de $S^2$, sont homologues (et connexes). Ainsi, $C_+$ (resp. $C_-$) est orienté de sorte qu'avec l'orientation induite, la composante $P_{k_1}^{-1}(\star)$ est homologue à l'opposé de la composante $P_{k_2}^{-1}(N)$ (resp. la composante $P_{k_1}^{-1}(-\star)$ est homologue à l'opposé de la composante $P_{k_2}^{-1}(S)$). Nous supposons dorénavant que l'orientation induite par celle de $C_+$ et celle de $C_-$ sur $P_{k_1}^{-1}(\star)$ et $P_{k_1}^{-1}(-\star)$ co\"\i ncide avec l'orientation naturelle comme image réciproque de $P_{k_1}$ et que l'orientation de $P_{k_2}^{-1}(N)$ et $P_{k_2}^{-1}(S)$ est opposée à celle obtenue via $P_{k_2}$.

On utilise l'anneau $\cA$, orienté, comme homologie de $C_+$ à zéro pour calculer l'enlacement entre $C_+(P_{k_1},P_{k_2})$ et $C_-(P_{k_1},P_{k_2})$ pour en déduire~: $$Enl(C_+(P_{k_1},P_{k_2}),C_-(P_{k_1},P_{k_2}))=k_1-k_2.$$

D'après le Théorème de Pontryagin (Théorème~\ref{t.pont}), $P_{k_1}$ et $P_{k_2}$ sont homotopes si et seulement si $k_1=k_2$ modulo $2p$, donc $f$ et $g$ sont homotopes si et seulement si $Enl(C_+(P_{k_1},P_{k_2}),C_-(P_{k_1},P_{k_2}))=Enl(C_+(f,g),C_-(f,g))=0$ modulo $2p$. Ce qui achève la preuve du lemme.
\end{proof}

Les lemmes~\ref{l.obst1geom2} et~\ref{l.obst2} permettent de prouver la Proposition~\ref{p.ns}.

\subsection{Exemple}\label{s.ex} Nous adaptons l'exemple de Pontryagin~\cite[\S4]{po} pour montrer l'utilité de notre critère. On considère les applications de $S^2\times S^1$ dans $S^2$\,; le Théorème~\ref{t.pont} assure l'existence de deux applications distinctes à homotopie près ayant le générateur de $H_1(S^2\times S^1)$ comme classe caractéristique.

\begin{figure}[htb]
\psfrag{x}{$x=\varphi_1(x)$}\psfrag{y}{$y$}\psfrag{z}{$z$}\psfrag{fy}{$\varphi_1(y)$}\psfrag{fz}{$\varphi_1(z)$}
\centerline{\includegraphics{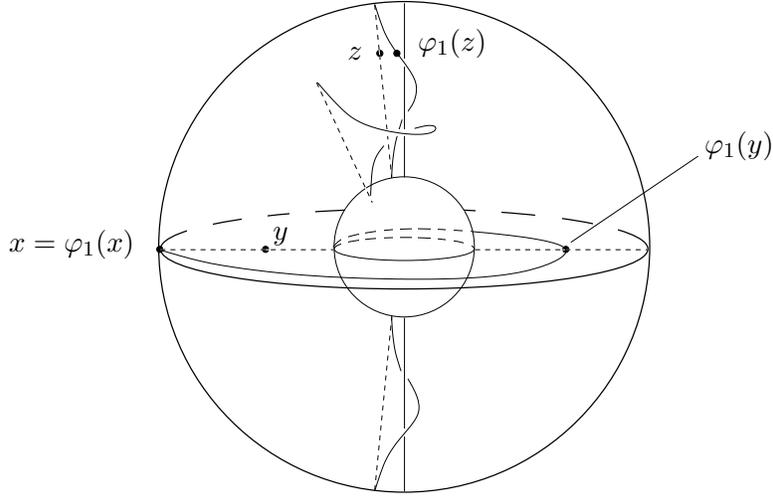}}
\caption{\label{f.dehn} Action de $\varphi_1$ sur $S^2\times S^1$}
\end{figure}

Une de ces applications est donnée par $\pi\colon S^2\times S^1 \to S^2$ la projection sur le premier facteur. Soient $p$ et $q$ deux points diametralement opposés de $S^2$, sphère unité de $\RR^3$, on note $\phi_\alpha$ la rotation d'angle $\alpha$ et d'axe $pq$. Dans les coordonnées polaires usuelles de l'espace $(r,\omega)$ ($r\in[0,\infty[$ et $\omega\in S^2$), on définit la suite de difféomorphismes $\varphi_n\colon \RR^3 \to \RR^3$ par $\varphi_n(r,\omega)=(r,\phi_{nr}(\omega))$. On note encore $\varphi_n$ la factorisation de $\varphi_n$ comme difféomorphisme de $S^2\times S^1$ (cf. figure~\ref{f.dehn}) et on remarque que $\varphi_n$ est homotope à $\varphi_m$ si et seulement si $n=m \mod 2$ ($\pi_1(SO(3))=\ZZ_2$ et $\varphi_1$ décrit un générateur de ce groupe). Ainsi, $\tilde f_n=\pi\circ\varphi_n\colon S^2\times S^1 \to S^2$ est homotope à $\pi$ si et seulement si $n$ est pair\,; $\pi$ et $\tilde f_1$ sont donc les deux seules applications, à homotopie près, ayant le générateur de $H_1(M)$ comme classe caractéristique.

Les applications $\tilde f_n$ et $\pi$ sont égales le long de deux cercles correspondants à l'axe de rotation et le long d'une sphère. Cette situation n'est pas générique, une petite perturbation de $\tilde f_n$ en $f_n=\phi_\varepsilon\circ\tilde f_n$ permet d'éliminer cette sphère.   

Les applications $f_n$ et $\pi$ sont opposées le long de $n$ cercles parallèles dans le plan $z=0$. Le calcul de l'orientation nous permet de conclure que $Enl(C_+(f_n,\pi),C_-(f_n,\pi))=n$ (au signe près). En appliquant le Lemme~\ref{l.obst2}, on retrouve le résultat~: $f_n$ est homotope à $\pi$ si et seulement si $n$ est pair.

\section{Application aux champs de vecteurs non-singuliers}\label{s.vecteurs}

Si $M$ est une variété de dimension trois compacte, sans bord et orientable, son fibré tangent $TM$ est trivialisable (voir~\cite[Problem 12.B]{mist} par exemple). En choisissant une trivialisation $\tau\colon TM\to M\times \RR^3$ on identifie les classes d'homotopie des champs de vecteurs non-singuliers sur $M$ (homotopie à travers les champs non-singuliers) avec les classes d'homotopie des applications de $M$ dans $S^2$. Cette remarque pose deux problèmes~: a priori, l'identification dépend de la trivialisation choisie\,; de plus, la construction explicite de trivialisations du fibré tangent est en général un problème difficile. Les résultats de la section précédente vont nous permettre de contourner en partie ces deux problèmes. 

Soit $X$ un champ de vecteurs non-singulier sur $M$ et $\tau$ une trivialisation du fibré tangent\,; on note $X_\tau\colon M \to S^2$ l'application induite et $c_{X_\tau}$ la classe d'homologie de l'image réciproque d'une valeur régulière de $X_\tau$. 

Rappelons que la {\em classe d'Euler} d'un champ de plans co-orienté sur une variété de dimension trois est l'obstruction (dans $H^2(M,\ZZ)$) à compléter une trivialisation de ce champ de plans du 1-squelette d'une triangularisation de la variété au 2-squelette. Nous appelons aussi classe d'Euler son dual de Poincaré (dans $H_1(M,\ZZ)$) et nous définissons la classe d'Euler d'un champ de vecteurs non-singulier, notée $\cE(X)$, comme étant la classe d'Euler du champ de plans orthogonal au champ de vecteurs, $X^\bot$, pour une métrique riemanienne quelconque. 

Un résultat classique de topologie algébrique nous permet de calculer la classe d'Euler de la façon suivante~: une section générique d'un champ de plans va couper la section nulle le long d'un entrelac orienté, la classe d'homologie de cet entrelac ne dépend pas du choix de la section générique, c'est la classe d'Euler du champ de plans.

Le résultat suivant est connu, afin d'\^etre complet, nous incluons la preuve de~\cite[Lemma 6.1.4]{bepe}.
\begin{lemme}
Pour toute trivialisation $\tau$ et tout champ de vecteurs $X$, on a la formule
$$\cE(X)=2 c_{X_\tau}.$$
\end{lemme}
\begin{proof}
Construisons une section générique de $X^\bot$, pour cela on choisit $y$ sur $S^2$ de sorte que $y$ et $-y$ soient des valeurs régulières de $X_\tau$. On définit $s\colon M \to X^\bot$ par $s(x)=\tau^{-1}(x,X_\tau(x)\wedge y)$ ($\wedge$ représente ici le produit vectoriel usuel de $\RR^3$). Cette section s'annule exactement quand $X_\tau(x)$ est colinéaire à $y$, c'est-à-dire sur l'ensemble $X_\tau^{-1}(y)\cup X_\tau^{-1}(-y)$. L'homologie de l'intersection de cette section avec la section nulle est donc égale à $2c_{X_\tau}$.
\end{proof}

\begin{remarque}\begin{itemize}
\item La classe d'Euler est invariante à homotopie près et indépendante du choix de la trivialisation du fibré tangent\,; ce n'est donc pas le cas de $c_{X_\tau}$ qui peut dépendre du choix de la trivialisation si $H_1(M,\ZZ)$ contient des éléments d'ordre 2.
\item Le diviseur maximal de $c_{X_\tau}$ ne dépend pas du choix de la trivialisation car c'est la moitié du diviseur maximal de $\cE(X)$.
\end{itemize}
\end{remarque}

\begin{definition}
Un champ de vecteurs $X$ est {\em complétable} si $\cE(X)=0$.
\end{definition}

Pour tout champ de vecteurs complétable, il existe une trivialisation $\tau_X$ telle que $X_{\tau_X}$ soit une application constante.

On retrouve facilement un joli résultat de Gompf~\cite[Corollary 4.10]{go}.

\begin{lemme}
Un champ de vecteurs $X$ est complétable si et seulement si $X$ est homotope à $-X$.
\end{lemme}
\begin{proof}
La première implication est évidente. Si $X$ est homotope à $-X$, pour n'importe quelle trivialisation $\tau$, $c_{X_\tau}=-c_{-X_\tau}=c_{-X_\tau}$. Ainsi $2c_{-X_\tau}=0$ et $2c_{X_\tau}=0$\,; donc $\cE(X)=0$.
\end{proof}

Soient $X$ et $Y$ deux champs de vecteurs non-singuliers sur $M$, on définit $C_+(X,Y)=\{x\in M,\, X(x)=\lambda Y(x)\mbox{ avec }\lambda>0 \}$ et $C_-(X,Y)=\{x\in M,\, X(x)=\lambda Y(x)\mbox{ avec }\lambda<0 \}$. Comme conséquence des résultats de la section précédente, on obtient la proposition et le corollaire suivants qui donnent des outils indépendants du choix de la trivialisation~:
\begin{proposition}\label{p.champ}
Deux champs de vecteurs non-singuliers $X$ et $Y$ sur $M$ sont homotopes si et seulement si
\begin{itemize} 
\item $[C_-(X,Y)]_{H_1 (M)}=0$, ce qui entra\^\i ne l'égalité des diviseurs maximaux des classes d'Euler de $X$ et $Y$ (on note $p$ cet entier) et
\item $Enl(C_+(X,Y),C_-(X,Y))=0$ modulo $p$.
\end{itemize}
\end{proposition}

\begin{corollaire}\label{c.champ}
\'Etant donnés deux champs de vecteurs $X$ et $Y$ sur $M$, 
si $[C_-(X,Y)]_{H_1 (M)}=0$ et $Enl(C_+(X,Y),C_-(X,Y))=0$, $X$ et $Y$ sont homotopes.
\end{corollaire} 

\section{Champs de Morse-Smale sur les variétés de Seifert}

\subsection{Construction initiale}\label{s.initiale}
On considère une variété de Seifert $M^3$ de base une surface compacte $S$ et de projection $p\colon M \to S$, on note $x_i$ pour $i$ allant de 1 à $n$ les points de $S$ où se projettent les fibres singulières de $p$. Si $Y_0$ est un champ de vecteurs Morse-Smale sans orbite périodique sur $S$, tel que tout $x_i$ soit une singularité de $Y_0$ (de type puit ou source), on choisit une métrique riemannienne sur $M$ et on relève $Y_0$ en un unique champ de vecteurs lisse, $X_0$, orthogonal aux fibres de $p$. Le champ $X_0$ est nul le long de toutes les fibres se projettant sur des  zéros de $Y_0$. On note $X_1$ le champ de vecteurs unitaire, tangent aux fibres et on pose $X=X_0+X_1$. Le champ de vecteur $X$ est un champ de Morse-Smale non-singulier lisse sur $M$, dont la dynamique se projette via $p$ sur celle de $Y_0$. 

\begin{remarque}
Quitte à rajouter des extrémas locaux (en particulier aux points $x_i$), une petite perturbation du champ de gradient d'une fonction de Morse sur $S$ permet de construire un champ $Y_0$ comme ci-dessus (sans orbite périodique).
\end{remarque}

\subsection{Opération de Wada}

Un entrelac dans $M$ est~{\em indexé} si on attribue à chaque composante un indice, 0, 1 ou 2. L'entrelac des orbites périodiques d'un champ de Morse-Smale est naturellement indexé par la dimension de la variété instable du point fixe de l'application de premier retour sur un disque transverse. 

Dans~\cite{wa}, Wada donne une caractérisation des entrelacs indexés réalisables comme entrelacs d'orbites périodiques de champs de Morse-Smale non singulier sur la sphère $S^3$. Il définit 6 opérations sur les entrelacs indexés de la sphère et un générateur : l'entrelac de Hopf indexé par 0 et 2. 

\begin{figure}[htb]
\centerline{\includegraphics{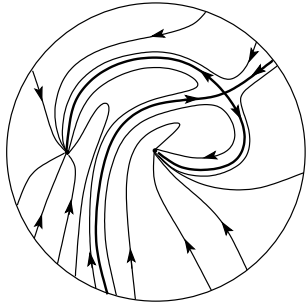}}
\caption{\label{f.Wada5.1} Une selle et deux attracteurs}
\end{figure}

Nous nous intéressons plus particulièrement ici à la cinquième opération de Wada qui consiste à remplacer un voisinage tubulaire d'une orbite périodique $K$, d'indice 0 ou 2, par un tore solide contenant trois orbites périodiques $K_1$, $K_2$ et $K_3$. L'\^ame du tore est l'une de ces orbites périodiques, $K_1$, son plongement dans $M$ à donc le m\^eme type de n\oe ud que $K$. Les orbites $K_2$ et $K_3$ sont des $n$-c\^ables parallèles de $K_1$. L'orbite $K_2$ est une selle et les indices de $K_1$ et $K_3$ sont 0 ou 2 mais au moins l'un des deux doit avoir l'indice de $K$. La suspension du difféomorphisme du disque décrit sur la figure~\ref{f.Wada5.1} permet d'obtenir une selle et un attracteur c\^ablant l'attracteur de départ (centre du disque).

Si le champ initial est de type Morse-Smale, il en est de m\^eme du champ obtenu en appliquant la cinquième opération de Wada.

\begin{remarque}
Le fibré normal de l'orbite $K$ n'est pas nécessairement trialisable de manière canonique (contrairement au cas de la sphère $S^3$), on ne peut donc pas parler de $(p,q)$-c\^ables autour de $K_1$ mais seulement de $n$-c\^ables à priori (c.f. \S\ref{s.enlace}).
\end{remarque}

\subsection{Orbites périodiques sur la surface de base}

\begin{figure}[htb]
\centerline{\includegraphics{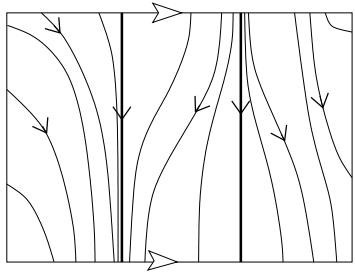}}
\caption{\label{f.tore} Dynamique sur les tores invariants}
\end{figure}

Si dans la construction initiale (\S~\ref{s.initiale}), le champ de Morse-Smale $Y_0$ admet des orbites périodiques sur $S$, le relevé $X$ n'est pas nécessairement Morse-Smale. Au-dessus de chaque orbite périodique, il laisse invariant un tore (ou une bouteille de Klein) attractif ou répulsif suivant la nature de l'orbite périodique de $Y_0$. On modifie alors la dynamique sur ce tore (resp. cette bouteille de Klein) de la façon suivante. On ajoute 2 orbites périodiques, de pente $(1,0)$ si la fibre de Seifert représente $(0,1)$ sur le tore. Pour la dynamique sur le tore, une de ces orbites est attractive, l'autre est répulsive (figure~\ref{f.tore}). On peut plonger ce tore dans le tore épaissi voisinage du tore invariant de départ de sorte que l'orbite répulsive du tore devient une selle pour $X$. 

La construction sur les bouteilles de Klein (cas où la surface de base n'est pas orientable) est similaire, on identifie les c\^otés verticaux de la figure~\ref{f.tore} par un homéomorphisme renversant l'orientation verticale.

L'un des points essentiels de la preuve du Théorème~\ref{t.nms} est la construction de champs de vecteurs $Y_0$ sur la surface de base ayant des orbites périodiques dont la classe d'homologie est prescrite. Si la surface de base est orientable, nous faisons les choix suivants.
 
\begin{lemme}\label{l.lecteur}
\begin{itemize}
\item Tout élément primitif de $H_1(\TT^2,\ZZ)$ peut \^etre représenté par une courbe fermée simple, plongée dans $\TT^2$.
\item Si $S=\#_{i=1}^{g}\TT^2$, pour tout élément $c$ de $H_1(S,\ZZ)$, il existe $\lambda_i$ pour $i$ allant de 1 à $g$ tels que $c=\sum_{i=1}^{g} \lambda_i[\gamma_i]$ où les $[\gamma_i]$ sont des éléments primitifs de $H_1(\TT^2,\ZZ)$ et les $\gamma_i$ sont deux à deux disjoints.
\item Soit $M$ une variété de Seifert au-dessus de $S=\#_{i=1}^{g}\TT^2$. Pour tout $c$ de $H_1(M,\ZZ)$, il existe $\lambda_i$ pour $i$ allant de 1 à $g$, $\lambda$ et $\alpha_j$ pour $j$ allant de 1 à $n$ tels que $$c=\sum_{i=1}^{g} \lambda_i[\gamma_i]+\lambda [F]+\sum_{j=1}^{n}\alpha_j[F_j]$$ où $F$ est une fibre générique, les $F_j$ sont les fibres singulières et les projections des $\gamma_i$ sur $S$ sont comme dans le point précédent. 
\end{itemize}
\end{lemme}

\begin{proof}
Les n\oe uds toriques $(p,q)$ avec $p$ et $q$ premiers entre eux représentent, par une courbe fermée simple, les éléments primitifs de $H_1(\TT^2,\ZZ)$. Le second point consiste simplement à effectuer la somme connexe des tores en dehors des $\gamma_i$.

Pour montrer le dernier point, nous remarquons que si l'on note $\tilde M$ la variété obtenue en enlevant un voisinage tubulaire de chaque fibre singulière et d'une fibre régulière, l'inclusion de $\tilde M$ dans $M$ induit un homomorphisme surjectif de  $H_1(\tilde M)$ dans $H_1(M)$. La fibration de $\tilde M$ est triviale, par un abus de notation, on pose $\gamma_i$ le relevé dans $\tilde M$ de $\gamma_i\subset \tilde F$ (où $\tilde F$ est la surface $F$ privée de $n+1$ points).
\end{proof}

Nous laissons au lecteur le soin d'adapter ce lemme au cas où la surface de base est non-orientable.

\begin{remarque}\label{r.preborne}
Le nombre de termes dans l'expression de $c$ au dernier point du Lemme~\ref{l.lecteur} ne dépend que de la variété de Seifert, il est inférieur ou égal à $g+n+1$.
\end{remarque}

\'Etant donné un entrelac sur la surface $S$, il est facile de construire un champ de Morse-Smale $Y_0$ sur $S$ ayant cet entrelac inclu dans l'ensemble de ses orbites périodiques attractives. En conséquence, on montre (en laissant encore au lecteur le soin de traiter le cas des surfaces non-orientables)~:

\begin{lemme}\label{l.realise}
Pour $M$ une variété de Seifert, tout élément de $H_1(M,\ZZ)$ est réalisable par la classe d'homologie d'un entrelac d'orbites périodiques attractives d'un champ de Morse-Smale non-singulier sur $M$.
\end{lemme}

\begin{proof}
\'Ecrivons $c\in H_1(M,\ZZ)$ comme la somme $c=\sum_{i=1}^{g} \lambda_i[\gamma_i]+\lambda [F]+\sum_{j=1}^{n}\alpha_j[F_j]$\,; on construit un champ de vecteurs $Y_0$ sur la surface de base $S$ ayant les projections des $\gamma_i$ (dans $S\setminus\{x_i\}$) comme orbites périodiques, $n+1$ singularités de type puit aux points $x_i$ et $x\in S\setminus\{x_i\}\cup\gamma_i$.

On relève $Y_0$ en un champ de Morse-Smale non-singulier $\tilde X$ comme précédemment. On applique la cinquième opération de Wada aux orbites attractives qui correspondent aux $\gamma_i$ avec le coefficient $\lambda_i$ (i.e. on obtient des $\lambda_i$-c\^ables de $\gamma_i$), aux orbites qui correspondent aux $x_j$ avec le coefficient $\alpha_j$ et à l'orbite se projettant sur $x$, avec le coefficient $\lambda$. On obtient ainsi un champ non-singulier $X$ qui est Morse-Smale et qui admet un entrelac d'orbites périodiques attractives dont l'homologie est $c$.
\end{proof}

\begin{remarque}\label{r.borne}
D'après la remarque~\ref{r.preborne} et la preuve précédente, on déduit qu'il existe un nombre $n(M)$ tel que tout élément de $H_1(M,\ZZ)$ est réalisable par la classe d'homologie d'un entrelac d'orbites périodiques attractives d'un champ de Morse-Smale non-singulier sur $M$ ayant au plus $n(M)$ orbites périodiques.
\end{remarque}

\begin{remarque}
En utilisant le résultat de Morgan~\cite{mo}, on retrouve un résultat de Yano~\cite{ya0} affirmant que tout élément de $H_1(M,\ZZ)$, pour $M$ une variété de Seifert, est représenté par un entrelac graphé (i.e. dont le complément est une variété graphée).
\end{remarque}

\subsection{Preuve du Théorème~\ref{t.nms}}
La preuve est divisée en deux étapes, dans un premier temps, on considère $M$ une variété de dimension trois, nous ne supposons pas que $M$ est une variété de Seifert.

\begin{lemme}\label{l.orientation}
Si $X$ est un champ Morse-Smale non-singulier sur $M$ et $\gamma$ est une orbite périodique attractive de $X$, il existe un champ Morse-Smale non-singulier $Y$ co\"\i ncidant avec $X$ sur le complément d'un voisinage tubulaire de $\gamma$ et ayant $-\gamma$ comme orbite périodique attractive.
\end{lemme}

\begin{figure}[htb]
\psfrag{g}{$\gamma$}\psfrag{-g}{$-\gamma$}
\centerline{\includegraphics{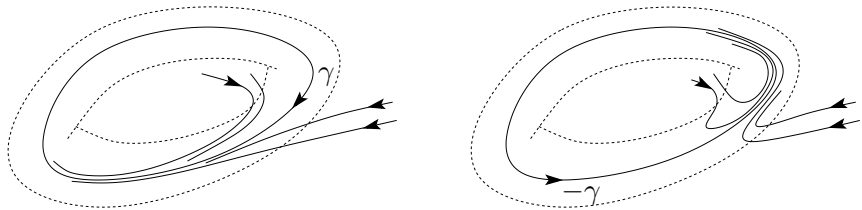}}
\caption{\label{f.orientation} Changement de l'orientation d'une orbite attractive}
\end{figure}

Le Lemme~\ref{l.obst1geom} nous permet de montrer que, pour $\tau$ une trivialisation du fibré tangent de $M$, les champs $X$ et $Y$ ci-dessus vérifient~: $$C_-(X,Y)=\gamma\mbox{ et ainsi }c_{X_\tau}=c_{Y_\tau}+[\gamma].$$

Rappelons que les champs $X$ tels que $[C_-(X_0,X)]=0$ sont ceux qui ont m\^eme <<demi-classe d'Euler>> que $X_0$.

\begin{proposition}\label{p.enlmorse}
Si $X_0$ est Morse-Smale, tout champ de vecteurs $X$ tel que $[C_-(X_0,X)]=0$ est homotope à un champ de vecteurs Morse-Smale.
\end{proposition}
\begin{proof}
La preuve consiste à construire un champ Morse-Smale dans chaque classe d'homotopie de champs $X$ telle que $[C_-(X_0,X)]=0$.

Soit $\gamma_0$ une orbite périodique attractive de $X_0$, on applique la cinquième opération de Wada à $\gamma_0$ pour obtenir $X_1$ ayant une orbite périodique attractive $\gamma_1$ et une orbite selle de plus que $X_0$ qui sont des $1$-c\^ables de $\gamma_0$. Les orbites $\gamma_0$ et $\gamma_1$ sont donc homotopes, on note $\cA$ l'anneau (que l'on choisit plongé dans $M$) bordant $\gamma_0\cup\gamma_1$. Les champs de vecteurs $X_0$ et $X_1$ sont homotopes.

Nous allons utiliser encore la cinquième opération de Wada sur $\gamma_0$ et $\gamma_1$, pour obtenir un champ de vecteurs possédant deux nouvelles orbites attractives $\gamma_0'$ et $\gamma_1'$, c\^ables respectifs de $\gamma_0$ et $\gamma_1$. Cette fois, l'existence de l'anneau $\cA$ nous permet de définir les invariants $(p_0,q_0)$ et $(p_1,q_1)$ des opérations de c\^ablage. Nous choisissons ces invariants comme suit : $(p_0=1,q_0=0)$ et $(p_1=1,q_1=\lambda)$ pour $\lambda$ un entier non nul et on note $X_2^\lambda$ le champ obtenu. Encore une fois, $X_2^\lambda$ est homotope à $X_0$.

La dernière étape de la construction consiste à changer l'orientation de $\gamma_0$ et $\gamma_1'$ en utilisant le Lemme~\ref{l.orientation} pour obtenir un champ $X_3^\lambda$. Une petite perturbation à l'extérieur d'un voisinage tubulaire de $\gamma_0$, $\gamma_1$, $\gamma_0'$ et $\gamma_1'$ permet d'obtenir: $$C_-(X_0,X_3^\lambda)=\gamma_0\cup\gamma_1'\mbox{ et}$$
$$C_+(X_0,X_3^\lambda)=\gamma_1\cup\gamma_0'.$$
Les champs de vecteurs Morse-Smale étant structurellement stables (c.f.~\cite{pasm}) $X_3^\lambda$ est Morse-Smale malgré la perturbation.

On a alors $$[C_-(X_0,X_3^\lambda)]=0\mbox{ pour tout $\lambda$; et}$$
$$Enl(C_+(X_0,X_3^\lambda),C_-(X_0,X_3^\lambda))=\lambda.$$

La Proposition~\ref{p.champ.ns} permet de conclure que l'on a effectivement construit au moins un champ Morse-Smale dans chaque classe d'homotopie ayant m\^eme <<demi-classe d'Euler>> que $X_0$.
\end{proof}

Pour la seconde étape de la preuve, nous considérons $M$ une variété de Seifert et nous construisons un champ de vecteurs Morse-Smale non-singulier ayant une <<demi-classe d'Euler>> donnée.  

\'Etant donnée $\tau$ une trivialisation du fibré tangent de $M$, on a~:

\begin{lemme}\label{l.euler}
Pour chaque $c$ dans $H_1(M,\ZZ)$, il existe un champ de Morse-Smale $X$, avec moins de $n(M)$ orbites périodiques, tel que $c_{X_\tau}=c$.
\end{lemme} 

\begin{proof}
Soit $Y_0$ un champ de Morse-Smale sur la surface $S$, sans orbite périodique et singulier aux points où les fibres singulières se projettent. On note $X_0$ le champ de vecteurs sur $M$ associé à $Y_0$ par la construction initiale. Le champ $X_0$ est homotope au champ de vecteurs tangent aux fibres de $M$, on note $e$ sa classe caractéristique pour la trivialisation $\tau$.

Soit $c$ un élément de $H_1(M,\ZZ)$, on note $c'=e-c$ et on considère $X'$ un champ de vecteurs Morse-Smale possédant un entrelac $\Gamma$ d'orbites périodiques attractives dont l'homologie est $c'$. Un tel champ, avec moins de $n(M)$ orbites périodiques, existe par le Lemme~\ref{l.realise} et la remarque~\ref{r.borne}.

On applique alors le Lemme~\ref{l.orientation} à $\Gamma$ pour obtenir un champ de vecteurs $X$ tel que $C_-(X',X)=\Gamma$. Le Lemme~\ref{l.obst1geom} permet d'obtenir $c_{X_\tau}=e-[\Gamma]=e-c'=c$.
\end{proof}

Le Théorème~\ref{t.nms} est une conséquence du Lemme~\ref{l.euler} et de la Proposition~\ref{p.enlmorse}.

\subsection{Discussion}\label{s.discussion}
On peut diviser les champs Morse-Smale non-singuliers d'une variété de dimension trois en trois catégories:
\begin{enumerate}
\item les champs dont aucune orbite périodique n'est homologue à zéro\,;
\item les champs ayant des orbites périodiques homologues à zéro mais vérifiant la {\em propriété d'enlacement} de Goodman (\cite{go1}, voir aussi~\cite{ya2}) : dès qu'une orbite périodique borde un disque, celui-ci rencontre une autre orbite périodique\,;
\item les champs ne vérifiant pas la propriété d'enlacement.
\end{enumerate}

Dans chaque classe d'homotopie de champs de vecteurs sur une variété de Seifert, nous avons construit une champ de Morse-Smale appartenant à la première catégorie. En particulier, d'après~\cite{go1,ya2}, ces champs de vecteurs sont transverses à un feuilletage.

{\bf Question~:} Peut-on classifier {\em à la Wada} les entrelacs indexés, sans composante homologue à zéro, d'une variété de Seifert réalisables comme entrelac d'orbites périodiques d'un Morse-Smale ? 

{\bf Question~:} Les autres entrelacs réalisables sont-ils obtenus par somme connexe avec un entrelac de Wada de $S^3$ ?

\bibliographystyle{halpha}
\bibliography{../champ}

\def\cprime{$'$}
\begin{thebibliography}{CMAV97}

\bibitem[Asi75]{as1}
Daniel Asimov.
\newblock Homotopy of non-singular vector fields to structurally stable ones.
\newblock {\em Ann. of Math. (2)}, 102(1):55--65, 1975.

\bibitem[BP97]{bepe}
Riccardo Benedetti and Carlo Petronio.
\newblock {\em Branched standard spines of {$3$}-manifolds}, volume 1653 of
  {\em Lecture Notes in Mathematics}.
\newblock Springer-Verlag, Berlin, 1997.

\bibitem[CGH03]{cogiho}
Vincent Colin, Emmanuel Giroux, and Ko~Honda.
\newblock On the coarse classification of tight contact structures.
\newblock to appear in the 2001 Georgia International Topology Conference
  proceedings, 2003.

\bibitem[CMAN98]{camanu}
J.~Casasayas, J.~Martinez~Alfaro, and A.~Nunes.
\newblock Knots and links in integrable {H}amiltonian systems.
\newblock {\em J. Knot Theory Ramifications}, 7(2):123--153, 1998.

\bibitem[CMAV97]{camavi}
B.~Campos, J.~Mart{\'{\i}}nez~Alfaro, and P.~Vindel.
\newblock Bifurcations of links of periodic orbits in non-singular
  {M}orse-{S}male systems on {$S\sp 3$}.
\newblock {\em Nonlinearity}, 10(5):1339--1355, 1997.

\bibitem[Duf03]{du_qtds}
Emmanuel Dufraine.
\newblock About homotopy classes of non-singular vector fields on the
  three-sphere.
\newblock to appear in Qualitative Theory of Dynamical Systems, 2003,
  arXiv:math.DS/0303315.

\bibitem[ET98]{elth}
Yakov~M. Eliashberg and William~P. Thurston.
\newblock {\em Confoliations}, volume~13 of {\em University Lecture Series}.
\newblock American Mathematical Society, Providence, RI, 1998.

\bibitem[Fra78]{fr}
John Franks.
\newblock The periodic structure of nonsingular {M}orse-{S}male flows.
\newblock {\em Comment. Math. Helv.}, 53(2):279--294, 1978.

\bibitem[Gei03]{ge}
Hansj\"org Geiges.
\newblock {Contact geometry}, 2003, arXiv:math.SG/0307242.

\bibitem[Gom98]{go}
Robert~E. Gompf.
\newblock Handlebody construction of {S}tein surfaces.
\newblock {\em Ann. of Math. (2)}, 148(2):619--693, 1998.

\bibitem[Goo85]{go1}
Sue Goodman.
\newblock Vector fields with transverse foliations.
\newblock {\em Topology}, 24(3):333--340, 1985.

\bibitem[Hon98]{hon}
Ko~Honda.
\newblock Confoliations transverse to vector fields.
\newblock preprint, 1998.

\bibitem[Kup94]{ku}
Krystyna Kuperberg.
\newblock A smooth counterexample to the {S}eifert conjecture.
\newblock {\em Ann. of Math. (2)}, 140(3):723--732, 1994.

\bibitem[Kup96]{ku1}
Greg Kuperberg.
\newblock Noninvolutory {H}opf algebras and {$3$}-manifold invariants.
\newblock {\em Duke Math. J.}, 84(1):83--129, 1996.

\bibitem[Mac01]{mk}
R.~S. MacKay.
\newblock Complicated dynamics from simple topological hypotheses.
\newblock {\em R. Soc. Lond. Philos. Trans. Ser. A Math. Phys. Eng. Sci.},
  359(1784):1479--1496, 2001.
\newblock Topological methods in the physical sciences (London, 2000).

\bibitem[Mil97]{mi}
John~W. Milnor.
\newblock {\em Topology from the differentiable viewpoint}.
\newblock Princeton Landmarks in Mathematics. Princeton University Press,
  Princeton, NJ, 1997.
\newblock Based on notes by David W. Weaver, Revised reprint of the 1965
  original.

\bibitem[Mor79]{mo}
John~W. Morgan.
\newblock Nonsingular {M}orse-{S}male flows on {$3$}-dimensional manifolds.
\newblock {\em Topology}, 18(1):41--53, 1979.

\bibitem[MS74]{mist}
John~W. Milnor and James~D. Stasheff.
\newblock {\em Characteristic classes}.
\newblock Princeton University Press, Princeton, N. J., 1974.
\newblock Annals of Mathematics Studies, No. 76.

\bibitem[NR90]{neru}
Walter~D. Neumann and Lee Rudolph.
\newblock Difference index of vectorfields and the enhanced {M}ilnor number.
\newblock {\em Topology}, 29(1):83--100, 1990.

\bibitem[Pon41]{po}
L.~Pontrjagin.
\newblock A classification of mappings of the three-dimensional complex into
  the two-dimensional sphere.
\newblock {\em Rec. Math. [Mat. Sbornik] N. S.}, 9 (51):331--363, 1941.

\bibitem[Pon59]{po2}
L.~S. Pontryagin.
\newblock Smooth manifolds and their applications in homotopy theory.
\newblock In {\em American Mathematical Society Translations, Ser. 2, Vol. 11},
  pages 1--114. American Mathematical Society, Providence, R.I., 1959.

\bibitem[PS70]{pasm}
J.~Palis and S.~Smale.
\newblock Structural stability theorems.
\newblock In {\em Global Analysis (Proc. Sympos. Pure Math., Vol. XIV,
  Berkeley, Calif., 1968)}, pages 223--231. Amer. Math. Soc., Providence, R.I.,
  1970.

\bibitem[Wad89]{wa}
Masaaki Wada.
\newblock Closed orbits of nonsingular {M}orse-{S}male flows on {$S\sp 3$}.
\newblock {\em J. Math. Soc. Japan}, 41(3):405--413, 1989.

\bibitem[Wil77]{wi}
F.~Wesley Wilson, Jr.
\newblock Some examples of nonsingular {M}orse-{S}male vector fields on
  {$S\sp{3}$}.
\newblock {\em Ann. Inst. Fourier. (Grenoble)}, 27(2):vi, 145--159, 1977.

\bibitem[Yan85a]{ya0}
Koichi Yano.
\newblock Homology classes which are represented by graph links.
\newblock {\em Proc. Amer. Math. Soc.}, 93(4):741--746, 1985.

\bibitem[Yan85b]{ya}
Koichi Yano.
\newblock The homotopy class of nonsingular {M}orse-{S}male vector fields on
  {$3$}-manifolds.
\newblock {\em Invent. Math.}, 80(3):435--451, 1985.

\bibitem[Yan85c]{ya2}
Koichi Yano.
\newblock Nonsingular {M}orse-{S}male flows on {$3$}-manifolds which admit
  transverse foliations.
\newblock In {\em Foliations (Tokyo, 1983)}, volume~5 of {\em Adv. Stud. Pure
  Math.}, pages 341--358. North-Holland, Amsterdam, 1985.

\end{thebibliography}

\end{document}